\documentclass[10pt,leqno]{amsart}
\usepackage{amssymb}
\usepackage{amsmath,amssymb,color}
\numberwithin{equation}{section}

\newcommand{\la}{\lambda}
\newcommand{\va}{\varphi}
\newcommand{\ppp}{\partial}

\newcommand{\DDD}{\mathcal{D}}
\newcommand{\www}{\widehat}

\newcommand{\pppa}{\partial_t^{\alpha}}

\newcommand{\R}{\mathbb{R}}

\newcommand{\C}{\mathbb{C}} 
\newcommand{\N}{\mathbb{N}}

\newcommand{\ooo}{\overline}
\newcommand{\OOO}{\Omega}

\newcommand{\sumij}{\sum_{i,j=1}^d}

\newcommand{\hhalf}{\frac{1}{2}}

\allowdisplaybreaks

\title
[]
{
Uniqueness of solution to boundary value problems for time-fractional wave 
equations
}

\author{
$^{1}$ Paola Loreti, $^1$ Daniela Sforza and $^{2,3,4}$ M. Yamamoto }

\thanks{
$^1$ Dipartimento di Scienze di Base e Applicate per l'Ingegneria, 
Sapienza Universit\`a di Roma,\\
Via Antonio Scarpa 16, 00161, Roma, Italy 
e-mail: {\tt paola.loreti@uniroma1.it}\\
{\tt daniela.sforza@uniroma1.it}\\
$^2$ Graduate School of Mathematical Sciences, The University
of Tokyo, Komaba, Meguro, Tokyo 153-8914, Japan \\
$^3$ Honorary Member of Academy of Romanian Scientists, 
Ilfov, nr. 3, Bucuresti, Romania \\
$^4$ Correspondence member of Accademia Peloritana dei Pericolanti,\\
Palazzo Universit\`a, Piazza S. Pugliatti 1 98122 Messina Italy
e-mail: {\tt myama@ms.u-tokyo.ac.jp}
}

\date{}

\begin{document}
\maketitle
\begin{abstract}
We consider an initial boundary value problem 
in a bounded domain $\OOO$ over a time interval $(0, T)$
for a time-fractional 
wave equation where the order of the fractional time derivative is
between $1$ and $2$ and the spatial elliptic operator has 
time-independent coefficients and is not necessarily symmetric.
We prove that if for arbitrarily chosen subdomain 
$\omega\subset \OOO$ and $T>0$,
a solution to the problem vanishes in $\omega \times (0,T)$ ,
then $u=0$ in $\OOO \times (0, T)$.
The uniqueness does not require any geometric condition on $\omega$.
\end{abstract}

\maketitle

\section{Introduction}

Let $\OOO \subset \R^d$ be a bounded domain with smooth boundary
$\ppp\OOO$, and let $1 < \alpha < 2$. 

For $1 < \alpha < 2$, let $d_t^{\alpha}$ denote 
the classical Caputo derivative:
$d_t^{\alpha}v(t) = \frac{1}{\Gamma(2-\alpha)}
\int^t_0 (t-s)^{1-\alpha} \frac{d^2v}{ds^2}(s) ds$
if $v, \frac{dv}{ds}, \frac{d^2v}{ds^2} \in L^1(0,T)$. Here 
$\Gamma(s)$ denotes the gamma function.
By $\pppa$ we denote an extension of $d_t^{\alpha}$ such that 
$\DDD(\pppa) \subset H^{\alpha}(0,T)$ which is the 
Sobolev-Slobodecki space (e.g., Adams \cite{Ad}).
Later in Section 1, we define $\pppa$.
Henceforth $L^2(\OOO)$, $H^1_0(\OOO)$, etc. denote usual 
Lebesgue space and Sobolev spaces.

We consider
$$
\left\{ \begin{array}{rl}
& \pppa (u(x,t) - a(x) - b(x)t) = -Au(x,t), \quad x\in\OOO, \, 0<t<T, \\
& u(x,\cdot) - a(x) - b(x)t \in H_{\alpha}(0,T) \quad 
   \mbox{for almost all $x \in \OOO$}, \\
& u(\cdot,t) \in H^1_0(\OOO) \quad \mbox{for almost all $t \in (0,T)$},
\end{array}\right.
                                         \eqno{(1.1)}
$$
where
$$
-Av(x) := \sumij \ppp_i(a_{ij}(x)\ppp_jv(x)) + \sum_{j=1}^db_j(x)\ppp_jv(x)
+ c(x)v(x), \quad x\in \OOO.                   \eqno{(1.2)}
$$
Here we assume that $a_{ij} = a_{ji} \in C^1(\ooo{\OOO})$, $b_j, c
\in C(\ooo{\OOO})$ for $1\le i,j \le d$ and $-A$ is uniformly 
elliptic. 
By $H^{\alpha}(0,T)$, we denote the Sobolev-Slobodecki space
(e.g., Adams \cite{Ad}).  Next, for $1<\alpha<2$, we define a Banach space
$H_{\alpha}(0,T)$ as follows:
$$
H_{\alpha}(0,T) := 
\left\{ \begin{array}{rl}
&\{ v\in H^{\alpha}(0,T);\, v(0) = 0\}, \quad 1 < \alpha < \frac{3}{2}, \\
&\left\{ v\in H^{\frac{3}{2}}(0,T); \thinspace
\int^T_ 0 t^{-1}\left\vert \frac{dv}{dt}(t)\right\vert^2 dt < \infty\right\}, 
\quad \alpha = \frac{3}{2}, \\
&\left\{v \in H^{\alpha}(0,T);\, v(0) = \frac{dv}{dt}(0) = 0\right\}, 
\quad \frac{3}{2} < \alpha < 2
\end{array}\right.
$$
with the norms $\Vert v\Vert_{H_{\alpha}(0,T)} = \Vert v\Vert_
{H^{\alpha}(0,T)}$ for
$\alpha \ne \frac{3}{2}$ and 
$\Vert v\Vert_{H^{\frac{3}{2}}(0,T)}
= \left( \Vert v\Vert^2_{H^{\frac{3}{2}}(0,T)} 
+ \int^T_0 t^{-1} \left\vert \frac{dv}{dt}(t)\right\vert^2 dt \right)
^{\hhalf}$.
We define the Riemann-Liouville fractional integral operator 
$J^{\alpha}$ by 
$J^{\alpha}v(t) := \frac{1}{\Gamma(\alpha)}\int^t_0
(t-s)^{\alpha-1} v(s) ds$ for $v\in L^2(0,T)$.
Then we can prove (e.g., Yamamoto \cite{Ya23})
that $J^{\alpha}: L^2(0,T)\, \longrightarrow H_{\alpha}(0,T)$ is 
isomorphism for $1<\alpha<2$ and 
see also Gorenflo, Luchko and Yamamoto \cite{GLY},
Kubica, Ryszewska and Yamamoto \cite{KRY}.  We define
$\pppa := (J^{\alpha})^{-1}$ with the domain $\DDD(\pppa)
= H_{\alpha}(0,T)$, and we can verify (e.g., \cite{Ya23}) that 
$\pppa v = d_t^{\alpha}v$ if $v, \frac{dv}{dt}, \frac{d^2v}{dt^2}
\in L^1(0,T)$ and $v(0) = \frac{dv}{dt}(0) = 0$.  Thus we can interpret that $\pppa$ is an extension of 
$d_t^{\alpha}$, and moreover the formulation with $\pppa u$ admits us to 
prove that for $a \in H^1_0(\OOO)$ and $b\in L^2(\OOO)$, there exists 
a unique solution $u \in H^1(0,T;L^2(\OOO)) \cap 
L^{\infty}(0,T;H^1_0(\OOO))$ to (1.1) such that 
$$
u-a-bt \in H_{\alpha}(0,T;H^{-1}(\OOO))          \eqno{(1.3)}
$$
(e.g., Huang and Yamamoto \cite{HY}).  
Here by $H^{-1}(\OOO)$, we denote the dual of $H^1_0(\OOO)$, identifying
$L^2(\OOO)$ with itself.
Throughout this article, we consider
solutions to (1.1) within this class.

The main subject of this article is
\\
{\bf Uniqueness}.
{\it 
Let $\omega$ be an arbitrarily chosen subdomain of $\OOO$.
Does we conclude that $a = b = 0$ in $\OOO$ if the solution u to (1.1) 
satisfies $u = 0$ in $\omega \times (0, T)$?
}

Although we do not discuss, the uniqueness is closely related to 
the approximate controllability through the duality.

The main purpose is to establish an affirmative answer:
\\
{\bf Theorem 1.}
\\
{\it
Let $a \in H^1_0(\OOO)$ and $b \in L^2(\OOO)$.  If a solution 
$u$ to (1.1) satisfies $u=0$ in $\omega \times (0, T)$,
then $u=0$ in $\OOO\times (0,T)$.
}

We emphasize that we do not require any geometric condition 
on $\omega$, which is essentially different from the case of the 
wave equation $\alpha=2$.
\\
{\bf Remark.}
In this article, for initial boundary value problem, we adopt 
the formulation (1.1).  
The classical formulation is based on the equation
$d_t^{\alpha}u(x,t) = -Au(x,t)$.
However, especially for non-symmetric case (1.2) with $a\in H^1_0(\OOO)$ and
$b \in L^2(\OOO)$, it is not direct to verify that $\ppp_t^2u(x,\cdot)
\in L^1(0,T)$ for fixed $x\in \OOO$, which can justify the definition of 
$d_t^{\alpha}u$ with $1<\alpha<2$.
Thus we adopt formulation (1.1).
There are other possible approaches for the well-posedness, and here 
we are restricted to Zacher \cite{Za}.  In the case of a symmetric 
operator, see Loreti and Sforza \cite{LS}.

Furthermore we can prove
\\
{\bf Lemma 1.}
{\it
Let $u(x,t)$ be the solution to (1.1) with 
$a\in H^1_0(\OOO)$ and $b\in L^2(\OOO)$.  \\
(i) The solution $u(x,z)$ to (1.1) is holomorphically 
extended to Re $z>0$ for arbitrarily fixed $x\in \OOO$.
\\
(ii) There exists a constant such that 
$\Vert u(\cdot,t)\Vert_{L^2(\OOO)} \le Ce^{Ct}
(\Vert a\Vert_{H^1_0(\OOO)} + \Vert b\Vert_{L^2(\OOO)})$ for all
$a\in H^1_0(\OOO)$, $b\in L^2(\OOO)$ and all $t>0$.
}

The proof is similar, for example, to Li, Imanuvilov and Yamamoto \cite{LIY},
and we omit the details of the proof.

There are existing works on the uniqueness of solution by 
extra data, and here we refer only to most related articles, 
not intending any comprehensive references.
\begin{itemize}
\item
{\bf Case $0<\alpha<1$}. 
We can refer to many works and see for example, Sakamoto and Yamamoto 
\cite{SY}, Jiang, Li, Liu and Yamamoto \cite{JLLY}, 
Jiang, Li, Pauron and Yamamoto \cite{JLPY}.
\item
{\bf Case $1 < \alpha < 2$}. we can not find many works and refer to 
\cite{JLLY}.
\end{itemize}

Moreover except for \cite{JLLY} and \cite{JLPY},  
all the works discuss symmetric A, that is, 
$b_j=0$ for $1\le j \le d$ and rely on the
eigenfunction expansion, which requires the symmetry.
The non-symmetric term $\sum_{j=1}^d b_j\ppp_ju$ represents an advection,
and is meaningful physically.
For the non-symmetric A, we cannot utilize eigenfunction expansions. 
The work \cite{JLLY} treats non-symmetric case for $0<\alpha<1$, 
and transforms the problem to the corresponding parabolic equation 
through the Laplace transform, and the method in [1] can work 
for $1< \alpha < 2$, but the problem is reduced to a wave equation, and for
the uniqueness we have to assume some geometric condition on 
$\omega$. 
Here to the case $1<\alpha<2$, we modify the argument of \cite{JLPY}
and eliminate any geometric condition on $\omega$ in \cite{JLLY}.

We can prove the uniqueness with extra data
$\ppp_{\nu_A}u\vert_{\gamma\times (0,T)}$ with arbitrarily chosen 
subboundary $\gamma$ and $(u(\cdot,\mu_k)_{L^2(\OOO)})_{1\le k \le m}$
with weight functions $\mu_1, ..., \mu_m$ satisfying suitable condition.
We can argue similarly when we replace the zero Dirichlet boundary 
condition by other boundary conditions in (1.1).


This article is composed of three sections.
In Section 2, we show some lemmata for the proof of Theorem 1.
In Section 3, we complete the proof of Theorem 1.
\section{Preliminaries}

{\bf 2.1. Laplace transform of $\pppa u$}

We set 
$(Lw)(p) := \www w(p) : = \int^{\infty}_0 e^{-pt}w(t) dt$,
provided that the integral converges for Re $p > p_0$: some positive constant.

Then 
\\
{\bf Lemma 2.}
{\it
Let $w(x,\cdot) - a(x) - b(x)t \in H_{\alpha}(0,T)$ for any $T>0$ 
and almost all $x\in \OOO$ and  
$\vert w(x,\cdot) \vert e^{-pt} \in L^1(0,\infty)$ for 
$p > p_0$ and almost all $x\in \OOO$.
Then
$L(\pppa (u-a-bt))(x,p) = p^{\alpha}(Lu)(p) - p^{\alpha-1}a(x)
- p^{\alpha-2}b(x), \quad Re p>p_0,\, x\in \OOO$.
}

The formula itself is well-known as
$\www{d_t^{\alpha}}(x,p) = p^{\alpha}\www{u}(p) 
- p^{\alpha-1}u(x,0) - p^{\alpha-2}\ppp_tu(x,0)$ for
$p> \mbox{Re}\, p_0$ and $x\in \OOO$,
but this requires the justification of the traces 
$u(x,0)$ and $\ppp_tu(x,0)$ in a suitable sense, which 
is not straightforward.
\\
{\bf Proof.}
First we remark: let $w \in H_{\alpha}(0,T)$ for any $T>0$ and 
$\vert w\vert e^{-pt} \in L^1(0,\infty)$ for 
$p > p_0$.  Then
$$
L(\pppa w)(p) = p^{\alpha}\www {w}(p) \quad \mbox{for Re $p>p_0$}.
                                 \eqno{(2.1)}
$$
The proof is found in Kubica, Ryszewska and Yamamoto \cite{KRY} and
Yamamoto \cite{Ya23}.

Now we will complete the proof of Lemma 2.  Henceforth, by 
$\langle \cdot, \cdot\rangle$ we denote the duality pairing
$\langle u, \va\rangle = \,_{H^{-1}(\OOO)}\langle u, \va\rangle_{H^1_0(\OOO)}$ 
for $u \in H^{-1}(\OOO)$ and $\va \in H^1_0(\OOO)$.
We arbitrarily choose $\psi\in C^{\infty}_0(\OOO)$.  Setting
$u_{\psi}(t) := \langle u(\cdot,t) - a - bt,\, \, \psi\rangle$,
by (1.3) and Lemma 1 (ii), since 
$\pppa u_{\psi}(t) = \langle \pppa (u(\cdot,t)-a-bt),\, \psi\rangle$ for 
$t>0$, we see that $u_{\psi}\in H_{\alpha}(0,T)$ for any $T>0$ and 
we can find a constant $C>0$ such that 
$\vert u_{\psi}(t)\vert \le Ce^{Ct}$ for all $t>0$.
Therefore we can apply (2.1), so that
\begin{align*}
& L(\langle \pppa (u-a-bt)),\, \psi \rangle)(p)
= L(\pppa u_{\psi})(p)
= p^{\alpha}L(u_{\psi})(p)\\
=& p^{\alpha}L(\langle u(\cdot,t),\psi\rangle
- \langle a,\psi \rangle - \langle bt, \psi\rangle)(p)
= p^{\alpha}(L(\langle u(\cdot,t),\psi\rangle)(p)
- L(\langle a,\psi \rangle)(p) - L(t \langle b, \psi\rangle)(p))\\
=& p^{\alpha} \langle (Lu(\cdot,p) - p^{-1}q - p^{-2}b),\,
\psi\rangle \quad \mbox{for Re $p > p_0$.}
\end{align*}
Therefore
$\langle L(\pppa (u-a-bt)),\, \psi \rangle(p)
= \langle (p^{\alpha}Lu(\cdot,p) - p^{-1}a - p^{-2}b),\,
\psi\rangle$ for Re $p > p_0$ and all $\psi\in C^{\infty}_0(\OOO)$.
Thus the proof of Lemma 2 is complete.
$\blacksquare$

{\bf 2.2. Spectral properties of $A$}

By attaching the operator $A$ in (1.2) with the domain $\DDD(A)
= H^2(\OOO) \cap H^1_0(\OOO)$, we define an operator in $L^2(\OOO)$, which 
is denoted by the same notation $A$.
Then it is known (e.g., Agmon \cite{Ag}) that 
$\sigma(A)$, the spectrum of $A$ is composed entirely of 
eigenvalues: $\sigma(A) := \{ \la_n\}_{n\in \N} \subset \C$.
Taking sufficiently small circle $\gamma_n$ which is centered at 
$\la_n$ and does not include $\gamma_m$ for $m\ne n$, we can 
define
$$
P_na = \frac{1}{2\pi\sqrt{-1}}\int_{\gamma_n} (z-A)^{-1}a dz,
D_na = \frac{1}{2\pi\sqrt{-1}}\int_{\gamma_n} (z-\la_n)(z-A)^{-1}a dz
$$
for $a \in L^2(\OOO)$ (e.g., Kato \cite{Ka}).
Then by \cite{Ag} and \cite{Ka}, we see that 
$P_n, D_n$ are bounded linear operators on $P_nL^2(\OOO)$ to itself, and 
$d_n:= \mbox{dim}\, P_nL^2(\OOO) < \infty$, and
$$
P_n^2 = P_n, \quad D_n=(A-\la_n)P_n, \quad D_nP_n=P_nD_n, \quad
D_n^{d_n}P_n=0.                            \eqno{(2.2)}
$$
Then we can show
\\
{\bf Lemma 3.}
{\it
We assume that $D_n^{k_0}P_n\va = 0$ with some $k_0\in \N$.  Then
$(A-\la_n)D_n^{k_0-1}P_n\va = 0$.
}
\\
{\bf Proof.}
By (2.2), we have
$P_nD_n^{k_0-1}P_n\va = D_n^{k_0-1}P_n^2\va = D_n^{k_0-1}P_n\va$.
Then, using (2.2) again, we obtain
$(A-\la_n)D_n^{k_0-1}P_n\va = ((A-\la_n)P_n)D_n^{k_0-1}P_n\va 
= D_nD_n^{k_0-1}P_n\va = D_n^{k_0}P_n\va = 0$. $\blacksquare$
\section{Proof of Theorem 1}

We write $\www{u}(p) = \www{u}(\cdot,p)$, etc.

{\bf First Step.} By (1.1), applying Lemma 2, we have
$p^{\alpha}\www{u}(p) - p^{\alpha-1}a - p^{\alpha-2}b
= -A\www{u}(p)$,
that is,
$$
\www{u}(p) = (p^{\alpha}+A)^{-1}(p^{\alpha-1}a + p^{\alpha-2}b) \quad
\mbox{for Re $p > p_0$}.
$$
Lemma 1 (i) yields  
$\www{u}(p) = 0 \quad \mbox{in $\omega$ for Re $p > p_0$}$,
that is, $p^{\alpha-1}(p^{\alpha}+A)^{-1}a = -p^{\alpha-2}(p^{\alpha}+A)^{-1}b$
in $\omega$ for Re $p > p_0$.
Setting $\eta:= p^{\alpha}$, we have
$\eta^{\frac{1}{\alpha}}(A+\eta)^{-1}a = -(A+\eta)^{-1}b$ in $\omega$ 
for Re $\eta > \eta_0$,
where $\eta_0$ is some positive constant.

We choose $\psi \in C^{\infty}_0(\OOO)$ arbitrarily.  Then
$$
(-\eta)^{\frac{1}{\alpha}}((A-\eta)^{-1}a, \, \psi)_{L^2(\omega)}
= -((A-\eta)^{-1}b, \, \psi)_{L^2(\omega)}            \eqno{(3.1)}
$$
for Re $\eta < \eta_0$.

We set 
$f_{\eta}(\eta):= ((A-\eta)^{-1}a,\psi)_{L^2(\omega)}$ and 
$g_{\eta}(\eta):= ((A-\eta)^{-1}b,\psi)_{L^2(\omega)}$ for
$\eta \in \rho(A):= \C \setminus \sigma(A)$.
We will prove
$$
f_{\psi}(\eta) = g_{\psi}(\eta) = 0 \quad \mbox{for all $\eta \in \rho(A)$}.
                                   \eqno{(3.2)}
$$
We note that $f_{\psi}(\eta)$ and $g_{\psi}(\eta)$ are holomorphic in 
$\eta \in \rho(A)$.

First let the zero set $Z_{\psi}:= \{ \eta \in \rho(A);\, 
f_{\psi}(\eta) = 0\}$
has an accumulation point.  Then the holomorphicity yields (3.2) in this case.
Second let $Z_{\psi}$ have no accumulation points.  Then, by (3.1), noting
that $\rho(A) \subset \{ \mbox{Re}\, z < q_0\}$ with some constant
$q_0$, we see that 
$(-\eta)^{\frac{1}{\alpha}} = -\frac{g_{\psi}(\eta)}{f_{\psi}(\eta)}$
if Re $\eta < \eta_0$ and $\eta \not\in Z_{\psi}$.

Since $\sigma(A) \subset \{ \mbox{Re}\, z > q_0\}$,
$\sigma(A)$ is a countably infinite and discrete set in $\C$, and 
$-\frac{g_{\psi}(\eta)}{f_{\psi}(\eta)}$ is holomorphic 
in $\eta \in \rho(A) \setminus Z_{\psi}$, we see that 
$-\frac{g_{\psi}(\eta)}{f_{\psi}(\eta)}$ can be holomorphically continued 
to some open neighborhood of $\{ \eta > 0;\, \eta \in \rho(A)\}$, which 
means that $(-\eta)^{\frac{1}{\alpha}}$ is holomorphically continued 
to some open neighborhood of $\{ \eta > 0;\, \eta \in \rho(A)\}$, and in 
particular, the function $(-\eta)^{\frac{1}{\alpha}}$ has the same limit 
as $\eta$ approaches a positive number from the uppe and the lower upper 
complex
plane.  Since $0 < \frac{1}{\alpha} < 1$ by $1<\alpha<2$, this is impossible.
Therefore $Z_{\psi}$ has an accumulation point, so that 
(3.2) is seen also in this case.
Thus we complete the proof of (3.2).
$\blacksquare$

Hence
$((A-z)^{-1}a,\, \psi)_{L^2(\omega)} = ((A-z)^{-1}b,\, \psi)_{L^2(\omega)}
= 0$ for $z \in \rho(A)$.
Since $\psi \in C^{\infty}_0(\omega)$ is arbitrary, we reach 
$(A-z)^{-1}a = (A-z)^{-1}b  = 0$ in $\omega$ for all $z \in \rho(A)$.

{\bf Second Step.} We will prove:
\\
If $a \in L^2(\OOO)$ and $(A-z)^{-1}a = 0$ in $\omega$ for 
all $z \in \rho(A)$, then $a=0$ in $\OOO$.

Indeed, we directly see
$$
D_n^{\ell}P_na = \frac{1}{2\pi\sqrt{-1}}\int_{\gamma_n}
(z-\la_n)^{\ell} (z-A)^{-1} a dz = 0 \quad \mbox{in $\omega$ for 
$\ell \in \N$.}                                  \eqno{(3.3)}
$$
Since $D_n^{d_n}P_na = 0$ by (2.2), we have $(A-\la_n)(D_n^{d_n-1}P_na) = 0$
in $\OOO$ by Lemma 2.  Moreover $D_n^{d_n-1}P_na = 0$ in $\omega$ by (3.3).
The unique continuation for the elliptic operator $A-\la_n$ yields
$D_n^{d_n-1}P_na = 0$ in $\OOO$.

Next we see that $(A-\la_n)(D_n^{d_n-2}P_na) = D_n^{d_n-1}P_na = 0$ in 
$\OOO$.  By (3.3) with $\ell = d_n-2$, we have 
$D_n^{d_n-2}P_na = 0$ in $\omega$.  Therefore, we apply the unique
continuation to obtain $D_n^{d_n-2}P_na = 0$ in $\OOO$.
Continuing this argument, we reach $P_na = 0$ in $\OOO$ fo each $n\in \N$.

Since Span $\left( \bigcup_{n=1}^N P_nL^2(\OOO)\right)$ is dense in $L^2(\OOO)$
(e.g., Agmon \cite{Ag}), we see that $a=0$ in $\OOO$.
Similarly we can prove $b=0$.  Thus the proof of Theorem 1 is complete.
$\blacksquare$
%
%
%

{\bf Acknowledgements.}
This article was completed during the stay of the third authoe at 
Sapienza University of Rome in January-February 2023.
The work was supported by Grant-in-Aid for Scientific Research (A) 20H00117 
of Japan Society for the Promotion of Science.


\begin{thebibliography}{99} %

\bibitem{Ad}
R.A. Adams,
{\it Sobolev Spaces}, Academic Press, New York, 1975.

\bibitem{Ag}
S. Agmon, {\it Lectures on Elliptic Boundary Value Problems},
van Nostrand, Princeton, 1965. 

\bibitem{GLY}
R. Gorenflo, Y. Luchko and M. Yamamoto, 
Time-fractional diffusion equation in the fractional Sobolev
spaces, Fract. Calc. Appl. Anal. {\bf 18} (2015) 799-820.

\bibitem{HY}
X. Huang and M. Yamamoto,
Well-posedness of initial-boundary value problem for time-fractional 
diffusion-wave equation with time-dependent coefficients, preprint,
arXiv:2203.10448

\bibitem{Ka}
T. Kato, {\it Perturbation Theory for Linear Operators},
Springer-Verlag, Berlin, 1995.

\bibitem{JLLY}
D. Jiang, Z. Li, Y. Liu and M. Yamamoto, Weak unique continuation property 
and a related inverse
source problem for time-fractional diffusion-advection equations, 
Inverse Problems {\bf 33} (2017) 055013.

\bibitem{JLPY}
D. Jiang, Z. Li, M. Pauron and M. Yamamoto, Uniqueness for fractional 
nonsymmetric diffusion equations and an application to an inverse source 
problem, Math. Meth. Appl. Sci. 2022;1-13, doi:10.1002/mma.8644

\bibitem{KRY}
A. Kubica, K. Ryszewska and M. Yamamoto, {\it Time-Fractional Differential 
Equations: A Theoretical Introduction}, Springer-Verlag, Tokyo, 2020.


\bibitem{LIY}
Z. Li, O. Imanuvilov and M. Yamamoto.
Uniqueness in inverse boundary value problems
for fractional diffusion equations,
Inverse Problems {\bf 32} (2016) 015004



\bibitem{LS} P. Loreti and D. Sforza, Weak solutions for time-fractional 
evolution equations in Hilbert spaces, Fractal and Fractional  
\textbf{5} (2021) 138 
https://doi.org/10.3390/fractalfract5040138

\bibitem{SY}
K. Sakamoto and M. Yamamoto, Initial value/boundary value problems for 
fractional diffusion-wave
equations and applications to some inverse problems, J. Math. Anal. Appl. 
{\bf 382} (2011) 426-447.

\bibitem{Ya23}
M. Yamamoto, Fractional calculus and time-fractional differential equations: 
revisit and construction
of a theory, Mathematics, Special issue Fractional Integrals and Derivatives: 
True versus False,
https://www.mdpi.com/2227-7390/10/5/698

\bibitem{Za}
R. Zacher, Weak solutions of abstract evolutionary integro-differential 
equations in Hilbert spaces, Funkcialaj Ekvacioj {\bf 52} (2009) 1-18.

\end{thebibliography}
\end{document}